\documentclass{gtmon_a}
\pdfoutput=1
\usepackage{pinlabel}

\proceedingstitle{The interaction of finite-type and Gromov--Witten
invariants (BIRS 2003)}
\conferencestart{15 November 2003}
\conferenceend{20 November 2003}
\conferencename{The interaction of finite-type and Gromov--Witten
invariants}
\conferencelocation{Banff International Research Station, Banff, Alberta,
Canada}

\editor{David Auckly}
\givenname{David}
\surname{Auckly}

\editor{Jim Bryan}
\givenname{Jim}
\surname{Bryan}

\title{Conormal bundles, contact homology and knot invariants}

\author{Lenhard Ng}
\givenname{Lenhard}
\surname{Ng}
\address{Department of Mathematics\\
Stanford University\\\newline
Stanford\\
CA 94305\\USA}
\email{lng@math.stanford.edu}
\urladdr{http://alum.mit.edu/www/ng/}

\volumenumber{8}
\issuenumber{}
\publicationyear{2006}
\papernumber{6}
\startpage{129}
\endpage{144}

\doi{}
\MR{}
\Zbl{}


\received{23 December 2004}
\revised{2 February 2005}
\accepted{24 January 2005}
\published{22 April 2006}
\publishedonline{22 April 2006}
\proposed{}
\seconded{}
\corresponding{}
\editor{}
\version{}

\arxivreference{math.SG/0412330}

\begin{asciiabstract}
We summarize recent work on a combinatorial knot invariant called
knot contact homology. We also discuss the origins of this invariant
in symplectic topology, via holomorphic curves and a conormal bundle
naturally associated to the knot.
\end{asciiabstract}

\AtBeginDocument{\let\bar\wbar\let\tilde\wtilde\let\hat\what}

\makeatletter
\def\cnewtheorem#1[#2]#3{\newtheorem{#1}{#3}[section]
\expandafter\let\csname c@#1\endcsname\c@proposition}

\newtheorem{proposition}{Proposition}[section]
\cnewtheorem{theorem}[proposition]{Theorem}
\cnewtheorem{lemma}[proposition]{Lemma}
\cnewtheorem{corollary}[proposition]{Corollary}
\cnewtheorem{conjecture}[proposition]{Conjecture}
\cnewtheorem{question}[proposition]{Question}

\theoremstyle{difinition}
\cnewtheorem{definition}[proposition]{Definition}

\cnewtheorem{remark}[proposition]{Remark}
\cnewtheorem{example}[proposition]{Example}
\makeatother  %  move after \newtheorem block

\def\C{\mathbb{C}}
\def\L{\mathcal{L}}
\def\R{\mathbb{R}}
\def\O{\mathcal{O}}
\def\P{\mathbb{P}}
\def\Z{\mathbb{Z}}
\def\M{\mathcal{M}}
\def\A{\mathcal{A}}
\begin{document}

\begin{abstract}
We summarize recent work on a combinatorial knot invariant called
knot contact homology. We also discuss the origins of this invariant
in symplectic topology, via holomorphic curves and a conormal bundle
naturally associated to the knot.
\end{abstract}

\maketitle

%*********************************************************************
\section{Introduction}

String theory has provided a beautiful correspondence between
enumerative geometry and knot invariants; for details, see the survey
by Mari{\~n}o \cite{Mar} or other papers in the present volume. This
correspondence applies methods from physics and algebraic geometry to
a construction, described below, which is essentially symplectic.

To symplectic geometers, there is a natural way to study this same
construction by counting holomorphic curves. The symplectic approach
leads to a knot invariant which seems to be distinct from the
invariants from physics. This manuscript serves as a survey of work
by the author \cite{KBIOne,KBITwo,Fr} studying this invariant, which is
called knot contact homology.

The string-theoretic approach originated in a conjecture of 't Hooft
relating large--$N$ gauge theories and open string theories, and a
subsequent observation by Witten \cite{Wit} connecting Chern--Simons
theory and topological strings; it was pioneered in influential papers
of Gopakumar and Vafa \cite{GV} and, most importantly for our purposes,
Ooguri and Vafa \cite{OV}. We will not attempt to present a careful
description of the Ooguri--Vafa construction, but only an extremely brief
summary.

A knot $K$ in $S^3$ gives rise to a natural Lagrangian submanifold
$\L K$ in the symplectic manifold $T^*S^3$ which is the conormal
bundle to $K$:
\[
\L K = \{(x,\xi)\ |\ x\in K \text{ and } \langle \xi, v \rangle = 0
\ \text{for all}\ v\in T_x K\} \subset T^*S^3
\]
There is a way of deforming $T^*S^3$ into the manifold
$\O(-1)\oplus \O(-1)\rightarrow \P^1$, by collapsing the zero section
and then performing a ``small resolution'' of the resulting conifold
singularity. This deformation then gives rise to a Lagrangian submanifold
$\widetilde{\L} K \subset \O(-1)\oplus\O(-1)$ which agrees with $\L
K\subset T^* S^3$ away from the deformation. There is then a
correspondence between knot invariants of $K$ derived from
Chern--Simons theory, and relative Gromov--Witten enumerative
invariants of $\widetilde{\L}K$:
\[
\begin{array}{ccc}
T^*S^3 &
%\stackrel{\text{conifold
%transition}}{\longleftrightarrow}
\xrightarrow{\text{conifold transition}} &
\O(-1)\oplus\O(-1) \\
\L K & \longleftrightarrow & \widetilde{\L}K \\
\text{C--S knot invariants} & \longleftrightarrow & \text{G--W
holomorphic invariants}
\end{array}
\]
From the viewpoint of symplectic geometry, however, it is natural to
look at enumerative holomorphic-curve invariants on the
\textit{left} side of this correspondence, before the conifold
transition. This gives rise to the knot invariants which are the
subject of the present paper.

We now describe how the symplectic knot invariants arise. Consider a
general $n$--dimensional Riemannian manifold $M$ with a compact
submanifold $K\subset M$. The Ooguri--Vafa setup uses a knot $K$ in
$M=S^3$; for technical reasons, it is simpler in the symplectic
setup to use $M=\R^3$ instead. There is a canonical $1$--form
$\lambda$ on the total space of $T^*M$ which pairs the result of
the projections $T(T^*M) \rightarrow TM$ and $T(T^*M)
\rightarrow T^*M$, and the natural symplectic structure on
$T^*M$ is given by $\omega = -d\lambda$. (In local coordinates
$x_i$ on $M$ and $\xi^i$ in the corresponding cotangent directions,
$\lambda = \xi^i dx_i$ and $\omega = dx_i \wedge d\xi^i$.) It is
easy to check that the conormal bundle $\L K$ is Lagrangian in
$T^*M$; that is, $\L K$ has dimension $n$ and $\omega|_{\L K} =
0$.

\begin{figure}[ht!]
\labellist\small
\pinlabel $T^*M$ at 93 147
\pinlabel $M$ [r] at 91 81
\pinlabel $M$ [br] at 240 106
\pinlabel $K$ [r] at 272 82
\pinlabel $ST^*M$ [bl] at 293 152
\pinlabel $LK$ [t] at 285 46
\pinlabel $\L K$ [t] at 376 30
\endlabellist
\centerline{
\includegraphics[width=4.5in]{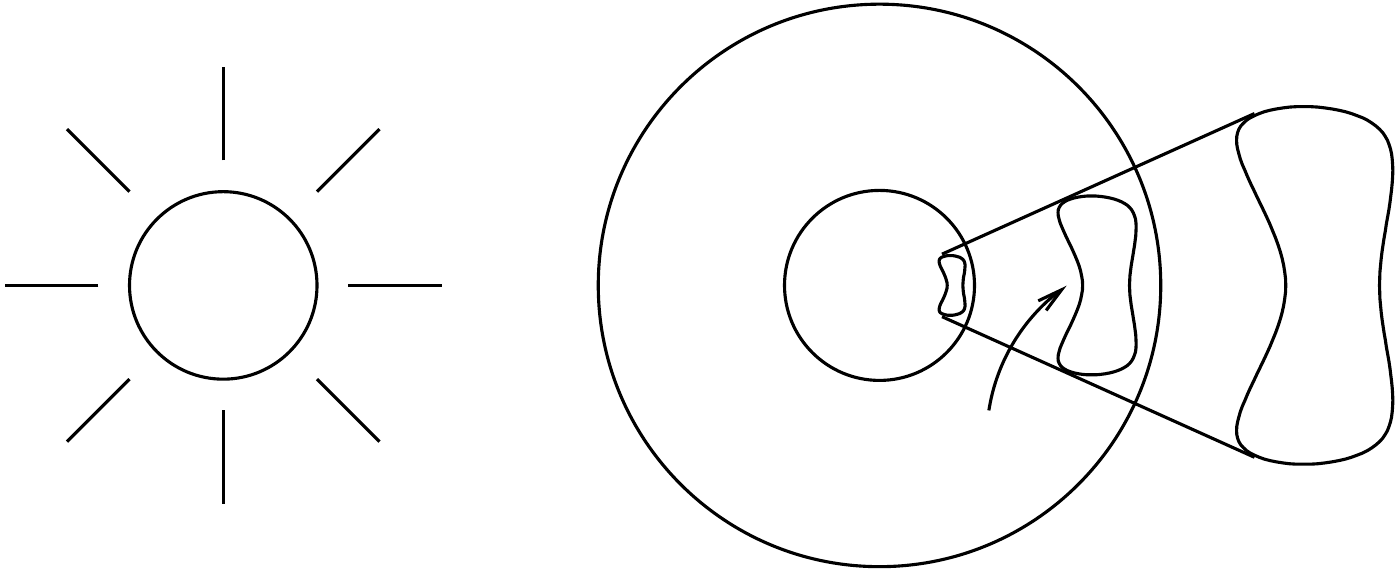}
}
\caption{
Schematic depiction of the cotangent bundle of $M$ (left) and the
conormals $\L K\subset T^*M$ and $LK\subset ST^*M$ (right). Note
that $\L K$ intersects $M$ in the knot $K$ and $ST^*M$ in the unit
conormal $LK$.
}
\label{fig-conormal}
\end{figure}
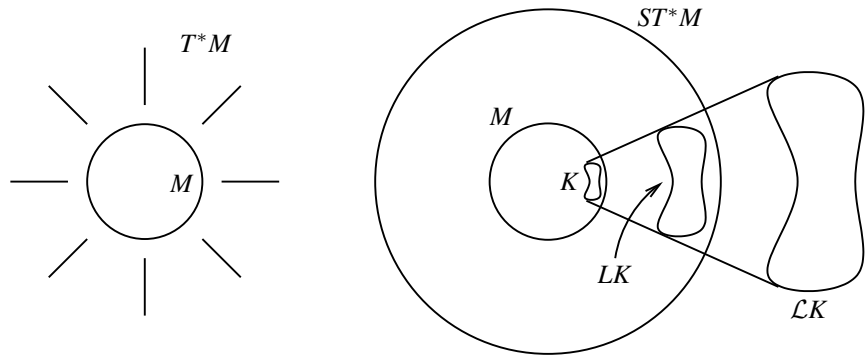

Geometrically, $\L K$ is somewhat awkward to use, because it is
noncompact. Instead, given a metric on $M$, we can work with the
\textit{unit} conormal bundle $LK$, which is the intersection of $\L
K$ with the cosphere bundle $ST^*M = \{(x,\xi)\in
T^*M\,|\,\|\xi\|=1\}$; see \fullref{fig-conormal}. The cosphere
bundle $ST^*M$ has a natural \textit{contact form}; recall that a
contact form, the odd-dimensional analogue of a symplectic form, is
a $1$--form $\alpha$ such that $\alpha\wedge d\alpha^{n-1}$ is
nowhere zero. In this case, the contact form is given by $\alpha =
\lambda|_{ST^*M}$. It should be noted that the \textit{contact
structure} on $ST^*M$ induced by $\alpha$, defined to be the
distribution $\ker\alpha$, is independent of the metric on $M$ and
hence depends only on the smooth structure on $M$. With respect to
this contact structure, the unit conormal $LK$ is
\textit{Legendrian}, meaning that $LK$ has dimension $n{-}1$ and
$\alpha|_{LK} = 0$.

Ambient isotopy of $K$ in $M$ leads to an isotopy of $LK$ in $ST^*M$
through Legendrian submanifolds; hence a Legendrian-isotopy invariant
of $LK$ yields a smooth-isotopy invariant of $K$, that is, a knot
invariant in the case where $K$ is a knot. A Legendrian-isotopy
invariant in contact topology is provided by \textit{Legendrian
contact homology}, Eliashberg \cite{Eli}, which in this case is the
simplest nontrivial approximation to the Symplectic Field Theory
\cite{EGH} of Eliashberg, Givental, and Hofer. We will define
Legendrian contact homology more carefully in \fullref{sec-LegCH},
but, roughly speaking, it counts holomorphic curves with boundary on a
Lagrangian cylinder over the Legendrian submanifold.

The knot invariant given by Legendrian contact homology, which we
call \emph{knot contact homology} and denote by $HC_*(K)$, takes the form
of a homology graded in degrees greater than or equal to $0$. We
describe several combinatorial forms for knot contact homology in
\fullref{sec-def}. For more details, see \cite{KBIOne,KBITwo,Fr}. It is
likely that knot contact homology has close ties to the ``string
topology'' introduced by Chas and Sullivan \cite{CS}.

Because of the similarity of the symplectic and string-theoretic
pictures, one might expect that knot contact homology incorporates
Chern--Simons knot invariants such as the Jones polynomial. Whether
this is true remains to be seen, but knot contact homology is at least
connected to ``classical'' knot invariants such as the Alexander
polynomial and the $A$--polynomial, Cooper--Culler--Gillet--Long
\cite{CCGLS}; through the latter, knot contact homology has ties to
$SL_2\C$--representations of the knot group. We discuss these
relations in \fullref{sec-apps}.

We remark that the symplectic approach produces isotopy invariants
of any submanifolds in any manifold, not just knots. What these
invariants might be for higher-dimensional ``knots'' is presently
unknown.

%*********************************************************************
%*********************************************************************
\subsection*{Acknowledgments}
I would like to thank Dave Auckly and Jim Bryan for organizing the
very illuminating BIRS workshop ``The interaction of finite type and
Gromov--Witten invariants,'' and the participants of the workshop
for helping me to streamline my presentation of the constructions
described in this manuscript. This work was supported by an American
Institute of Mathematics Five-Year Fellowship.

%*********************************************************************
%*********************************************************************
\section{Legendrian contact homology}
\label{sec-LegCH}

Let $V$ be a contact manifold with contact form $\alpha$, and let $L
\subset V$ be a Legendrian submanifold. The \textit{Reeb vector
field} $R_\alpha$ on $V$ is uniquely defined by the conditions
$\iota(R_\alpha) d\alpha = 0$, $\alpha(R_\alpha)=1$. A \textit{Reeb
orbit} in $V$ is a closed orbit of the flow under the Reeb vector
field, and a \textit{Reeb chord} of $L$ is a path along the flow of
the Reeb vector field which begins and ends on $L$. The
\textit{symplectization} of $V$ is the symplectic manifold $W =
V\times\R$ with symplectic form $\omega = d(e^t\alpha)$, where $t$
is the coordinate on $\R$ and $\alpha$ is induced from $V$. We can
give $W$ an almost complex structure by choosing any almost complex
structure on the distribution $\ker \alpha$, and further setting $J
\partial_t = R_\alpha$.

The Symplectic Field Theory of $V$ is a Floer-theoretic algebraic
structure derived from counting holomorphic curves in the
symplectization of $V$ which limit to Reeb orbits at $\pm \infty$.
As a more manageable approximation, we can study Eliashberg and
Hofer's contact homology \cite{Eli}, which counts genus--$0$
holomorphic curves with one end at $+\infty$ and an arbitrary number
of ends at $-\infty$. More precisely, the curves in question are
holomorphic maps from a multiply-punctured sphere to $V\times\R$
such that a neighborhood of one puncture tends to a cylinder over a
Reeb orbit as $t\to +\infty$, and neighborhoods of the other
punctures tend to cylinders over Reeb orbits as $t\to -\infty$.
%These curves arise naturally when considering the compactification
%of the moduli space of holomorphic genus--$0$ curves in $V\times\R$.
Out of the Reeb orbits and these holomorphic curves, one forms a
graded complex whose homology is called \emph{contact homology} and is an
invariant of the contact structure.

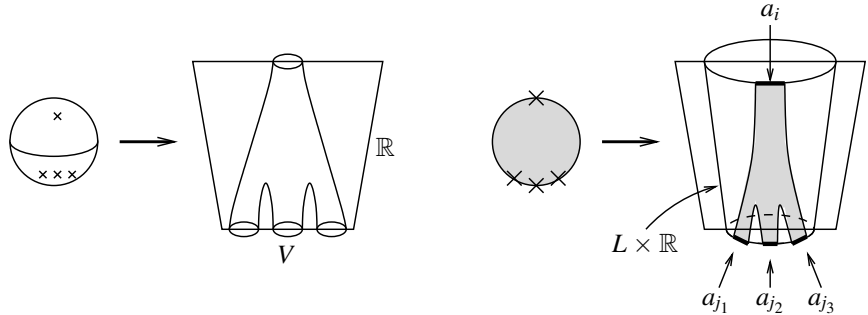
\begin{figure}[ht!]
\labellist\small
\pinlabel $V$ [t] at 128 24
\pinlabel $\R$ [l] at 166 66
\pinlabel $L\times\R$ [t] <0pt, 2pt> at 294 27
\pinlabel $a_i$ [b] at 352 121
\pinlabel $a_{j_1}$ [t] <0pt, 2pt> at 327 0
\pinlabel $a_{j_2}$ [t] <0pt, 2pt> at 352 0
\pinlabel $a_{j_3}$ [t] <0pt, 2pt> at 376 0
\endlabellist
\centerline{\includegraphics[width=4.5in]{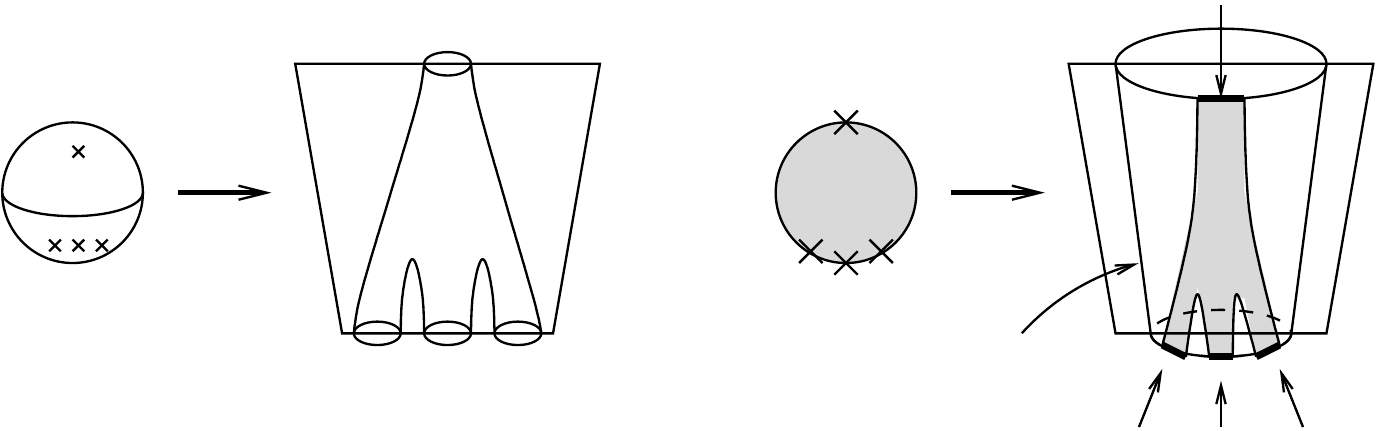}}\vspace{2mm}
\label{fig-CH}
\caption{
On the left, a holomorphic curve in $V\times\R$ contributing to the
contact homology of $V$. The punctures are mapped to the cylindrical
ends, which approach cylinders over Reeb orbits of $V$. On the
right, a holomorphic disk in $V\times\R$ contributing to the
Legendrian contact homology of $L$. Here the boundary of the disk,
away from the punctures, is mapped to $L\times\R$, and neighborhoods
of the punctures are mapped to strips approaching Reeb chords of
$L$. This disk contributes the term $a_{j_1}a_{j_2}a_{j_3}$ to $\partial
a_i$.
}
\end{figure}

Given a Legendrian submanifold $L\subset V$, one can construct a
relative version of contact homology called \textit{Legendrian
contact homology}. For simplicity, we will assume that the ambient
manifold $V$ has no closed Reeb orbits; this holds, for instance, in
the case $V = ST^*\R^3$. One then counts holomorphic disks which
limit to Reeb chords of $L$ at $t\to \pm \infty$, and whose boundary
lies on $L$.

More precisely, suppose that $L$ has a finite number of Reeb chords, which
we label $a_1,\dots,a_n$. The complex whose homology yields
Legendrian contact homology is the tensor algebra $\A$ freely
generated by $a_1,\dots,a_n$ over the group ring $\Z[H_1(L)]$; that
is, $\A$ is the module over $\Z[H_1(L)]$ generated by
(noncommutative) words in $a_1,\dots,a_n$, including the empty word.
This algebra can be given a grading by placing the base ring
$\Z[H_1(L)]$ in degree $0$, and assigning degrees to the generators
$a_i$ which are essentially the Conley--Zehnder indices of the Reeb
chords; see \cite{EES}. In order to define a differential on $\A$,
we need to study certain moduli spaces of holomorphic curves.

For each Reeb chord $a_i$, choose a ``capping path'' lying in $L$
which joins the endpoints of the chord. Let $D_k$ be a disk with
boundary, minus boundary punctures $x,y_1,\dots,y_k$ appearing in
order around the boundary, and give $D_k$ the complex structure
induced from the unit disk in $\C$. For $A\in H_1(L)$, we define the
moduli space $\M^A(a_i;a_{j_1},\dots,a_{j_k})$ to be the set of maps
$f\co D_k \to V\times\R$ such that:
\begin{itemize}
\item $f$ is holomorphic and proper and has finite energy with respect
  to $d\alpha$;
\item $f$ maps the boundary of $D_k$ (without the punctures) to
$L\times\R$;
\item $f$ maps a neighborhood of the puncture $x$ to a strip
approaching the strip $a_i\times\R$ as $t\to+\infty$;
\item $f$ maps a neighborhood of the puncture $y_i$ to a strip
approaching the strip $a_{j_i}\!\times\R$ as $t\to -\infty$;
\item the image of the boundary of $D_k$, made into a closed curve
by appending the appropriate capping paths, is in the homology class
$A$.
\end{itemize}
Note that there is an $\R$--action on
$\M^A(a_i;a_{j_1},\dots,a_{j_k})$ given by translation in the $\R$
direction. The moduli spaces that interest us are the ones that are
rigid modulo this $\R$--action.

Define the differential of $a_i$ as
\[
\partial a_i = \sum_{\dim \M^A(a_i;a_{j_1},\dots,a_{j_k}) = 1} \#
\bigl(\M^A(a_i;a_{j_1},\dots,a_{j_k})/\R\bigr) \,A \,a_{j_1}\dots a_{j_k},
\]
where $\#\bigl(\M^A(a_i;\!a_{j_1},\dots,a_{j_k})/\R\bigr)$ is the signed number
of points in the $0$--dimensional space $\M/\R$; extend $\partial$ to $\A$
via the Leibniz rule. Then $(\A,\partial)$ becomes a differential graded
algebra, usually abbreviated in the subject as a DGA.

\begin{theorem}
\label{thm-LegCH}
$\partial^2=0$, $\partial$ lowers degree by $1$, and the graded homology
$H_*(\A,\partial)$ is an invariant of the Legendrian isotopy class of $L$.
\end{theorem}

It is currently a bit of a misnomer to label this result, in full
generality, as a theorem, since the analytical details to prove it
are still being worked out. The foundational analysis has been
performed in several cases of interest, including $\R^3$ \cite{ENS},
$\R^{2n+1}$ for $n>1$ \cite{EES}, and jet spaces. On the other hand,
there is a standard technique  for
sidestepping these analytical issues by finding a purely
combinatorial form for the differential graded algebra $(\A,\partial)$ and
proving that this combinatorial algebra is invariant under
Legendrian isotopy. This technique dates back to Chekanov's pioneering
work \cite{Che} on Legendrian knots in standard contact $\R^3$, and is
the strategy used in \cite{KBIOne,KBITwo,Fr} to define knot contact homology.

In the case of interest here, $V$ is the contact manifold $ST^*\R^3
= \R^3\times S^2$, which can also be viewed as the $1$--jet space
$J^1(S^2)$, and $L$ is the Legendrian torus $LK$ given by the unit
conormal to the knot $K$ in $ST^*\R^3$. Given a framing and an
orientation of $K$, we obtain longitude and meridian classes
$\lambda,\mu\in H_1(LK)$ and hence an identification of
$\Z[H_1(LK)]$ with $\Z[\lambda^{\pm 1},\mu^{\pm 1}]$.

The Reeb vector field lies within the $\R^3$ fibers in $\R^3\times
S^2$; in the fiber lying over $\xi\in S^2$, it simply points in the
direction of the dual to $\xi$. It follows that Reeb chords of $LK$
correspond to ``binormal chords'': oriented line segments beginning
and ending on $K$ which are normal to $K$ at both endpoints. If we
parametrize $K$ by $S^1$, then the distance function in $\R^3$
between points on $K$ gives a map $d\co S^1\times S^1\to \R$ whose
nondiagonal critical points are binormal chords. One can then set up
the theory so that the degree of a Reeb chord in the contact
homology DGA is the Morse index of the corresponding critical point
of $d$. Hence the generators of the DGA all have degree $0$, $1$, or
$2$; it follows that the DGA is a nonnegatively graded algebra over
the ring $\Z[\lambda^{\pm 1},\mu^{\pm 1}]$.

\begin{figure}[ht!]
\labellist\small
\pinlabel $\R^3$ [tl] at 73 147
\pinlabel $K$ [t] at 150 60
\pinlabel $T^*\R^3$ at 169 13
\pinlabel $\L K$ [br] at 245 172
\pinlabel $a_i$ [l] at 312 129
\pinlabel $a_{j_1}$ [br] <0pt,-2pt> at 151 144
\pinlabel $a_{j_2}$ [r] at 146 122
\pinlabel $a_{j_3}$ [tr]  at 147 103
\endlabellist
\centerline{
\includegraphics[height=2.7in]{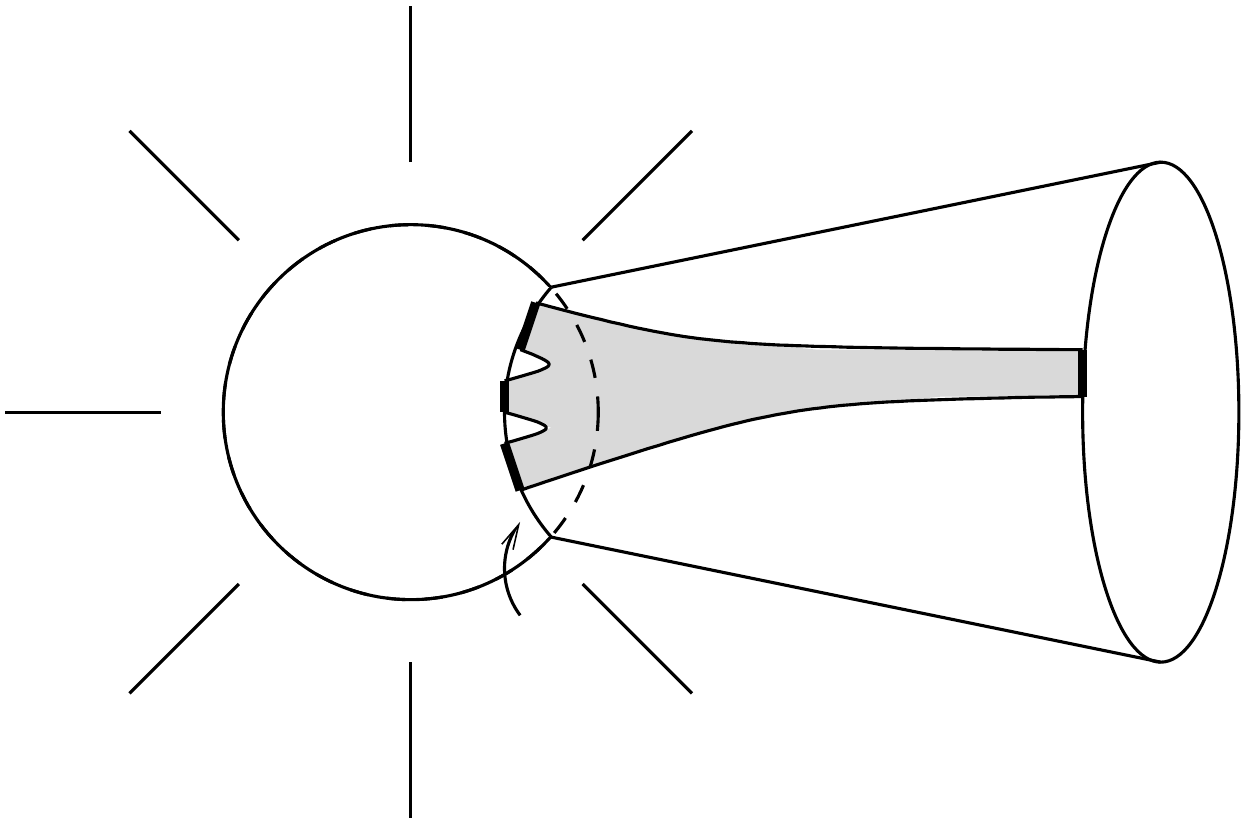}
}
\caption{
Holomorphic disk used in the computation of the Legendrian contact
homology of $LK$. The ambient space is $T^*\R^3$. In this picture,
the holomorphic disk has boundary on $\L K$, except for the strips
which limit to the binormal chords $a_{j_1},a_{j_2},a_{j_3}$ on the
zero section and the binormal chord $a_i$ at infinity.
}
\label{fig-LCH}
\end{figure}
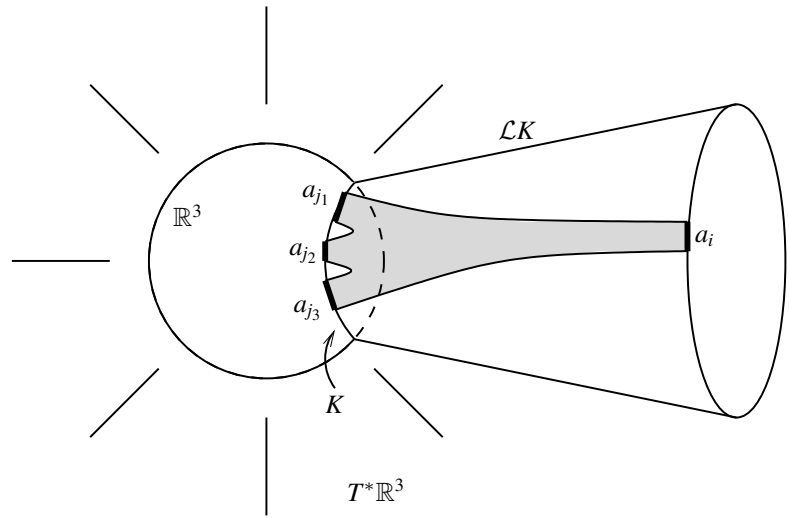

We can view Legendrian contact homology in $ST^*\R^3$ as counting
curves in $T^*\R^3$. The symplectization
$ST^*\R^3\times\R$ is symplectically diffeomorphic to
$T^*\R^3\setminus \R^3$, the cotangent bundle minus the zero
section, via the map sending $((x,\xi),t)$ to $(x,e^t\xi)$. The
$t\to+\infty$ end of the symplectization maps to the unbounded
``end'' of $T^*\R^3$, while the $t\to-\infty$ end maps to the zero
section. Under this identification, the holomorphic curves used to
calculate Legendrian contact homology are maps $f$ of
boundary-punctured disks to $T^*\R^3\setminus\R^3$ satisfying the
following properties:
\begin{itemize}
\item
$f$ is holomorphic, proper, and finite energy;
\item
$f$ sends the boundary of the disk to $\L K$;
\item
$f$ sends a neighborhood of one puncture to a strip which
approaches, on the unbounded end of $T^*\R^3$ (in the fiber
direction), a cylindrical strip of the form $BC \times \{e^t
\xi\ |\ t\gg 0\} \subset T^*\R^3$, where $BC$ is a binormal chord
and $\xi$ is a covector dual to the direction of $BC$;
\item
$f$ sends neighborhoods of the other punctures to neighborhoods of
binormal chords of $K$ in the zero section of $T^*\R^3$.
\end{itemize}
See \fullref{fig-LCH}.

We remark that the curves which produce Legendrian contact homology
in this context are somewhat different from the curves which would
be counted in relative Gromov--Witten theory; it is unknown what
relation our invariant has to the Chern--Simons type invariants
studied by Ooguri--Vafa et al. Also note that the almost complex
structure on $T^*\R^3\setminus \R^3$ inherited from the
symplectization, which we use to define our holomorphic curves,
cannot be extended to all of $T^*\R^3$.

This discussion of Legendrian contact homology should be viewed as
motivation for the knot contact homology invariant which will be
defined combinatorially in the following section. There is work in
progress to verify that the combinatorial DGA actually yields
Legendrian contact homology, by examining gradient flow trees on the
front of the Legendrian torus $LK$, and using the relation between
gradient trees and holomorphic curves described by Fukaya and Oh
\cite{FO}. As mentioned earlier, one can instead directly show that
the combinatorial theory gives a topological knot invariant, without
using its origin in symplectic geometry; this is the approach of
\cite{KBIOne,KBITwo,Fr}, and will be our approach in the next section.

%*********************************************************************
%*********************************************************************
\section{Knot contact homology: definition}
\label{sec-def}

Let $K$ be an oriented knot in $\R^3$; this has a canonical framing
given by any Seifert surface. As discussed in the previous section,
the Legendrian contact homology of $LK \subset ST^*\R^3$, which is a
knot invariant, is the homology of a complex defined in nonnegative
degree over the ring $\Z[\lambda^{\pm 1},\mu^{\pm 1}]$. We will now
give a combinatorial definition of this complex, which is called the
\textit{framed knot DGA}. (This complex may be different from the
one described in the previous section, but it has the same
homology.)

Fix a knot diagram for $K$ with $n$ crossings. The framed knot DGA
is the algebra $\A$ over $\Z[\lambda^{\pm 1},\mu^{\pm 1}]$ freely
generated by the following generators:
\begin{itemize}
\item $\{a_{ij}\}_{1\leq i,j\leq n,~i\neq j}$ of degree $0$;
\item $\{b_{\alpha i}\}_{1\leq \alpha,i\leq n}$ and
$\{c_{i\alpha}\}_{1\leq \alpha,i\leq n}$ of degree $1$;
\item $\{d_{\alpha \beta}\}_{1\leq \alpha,\beta\leq n}$ and
$\{e_\alpha\}_{1\leq \alpha\leq n}$ of degree $2$.
\end{itemize}
Here Greek subscripts represent crossings in the knot diagram,
numbered arbitrarily from $1$ to $n$, and Roman subscripts represent
components of the knot diagram, from undercrossing to undercrossing,
also numbered arbitrarily from $1$ to $n$.

To define the differential on $\A$, we first need some auxiliary
definitions. Each crossing $i$ involves three components of the knot
diagram, the overstrand $o_i$ and the understrands $l_i$ and $r_i$,
distinguished as being on the left and right sides as the overstrand
is traversed in the direction of the knot's orientation. Define
$\epsilon_1$ to be $\pm 1$ depending on the sign of crossing $1$,
where sign is defined in terms of the knot's orientation in the
usual way. (If $l_1$ follows $r_1$ when traversing the knot, then
$\epsilon_1=1$; if $r_1$ follows $l_1$, then $\epsilon_1=-1$.)

Let $\Psi^L,\Psi^R,\Psi^L_2,\Psi^R_1$ be the $n\times n$ matrices
defined as follows:
\begin{align*}
(\Psi^L)_{\alpha i} &= \begin{cases} \lambda^{-\epsilon_1} & \alpha=1, i=r_1 \\
1 & \alpha\neq 1, i=r_\alpha \\
\mu & i=l_\alpha \\
-a_{l_\alpha o_\alpha} & i=o_\alpha \\
0 & \text{otherwise} \end{cases} &\hspace{0.2in}
 (\Psi^R)_{i\alpha} &=
\begin{cases} \lambda^{\epsilon_1} \mu &
\alpha=1, i=r_1 \\
\mu & \alpha\neq 1, i=r_\alpha \\
1 & i=l_\alpha \\
-a_{o_\alpha l_\alpha} & i=o_\alpha \\
0 & \text{otherwise}
\end{cases}
\\
(\Psi^L_2)_{\alpha i} &= \begin{cases}
\mu & i=l_\alpha \\
-a_{l_\alpha o_\alpha} & i=o_\alpha \\
0 & \text{otherwise} \end{cases} &\hspace{0.2in}
 (\Psi^R_1)_{i\alpha} &=
\begin{cases} \lambda^{\epsilon_1} \mu &
\alpha=1, i=r_1 \\
\mu & \alpha\neq 1, i=r_\alpha \\
0 & \text{otherwise.}
\end{cases}
\end{align*}
Assemble generators of $\A$ into $n\times n$ matrices $A,B,C,D$ as
follows:
$$A_{ij} = \left\{ \begin{smallmatrix} 1+\mu & \text{if}\ i=j \\
                   a_{ij} & \text{if}\ i\neq j \end{smallmatrix} \right.;
  \quad B_{\alpha i} = b_{\alpha i};
  \quad C_{i\alpha} = c_{i\alpha};
  \quad\text{and}\ D_{\alpha\beta} = d_{\alpha\beta}.$$
Write $\partial A$ for the matrix whose entries are the
differentials of the corresponding entries of $A$, and similarly for
$\partial B$, $\partial C$, and $\partial D$. Then the differential on generators is
given by
\begin{align*}
\partial A &= 0 \\
\partial B &= \Psi^L \cdot A \\
\partial C &= A \cdot \Psi^R \\
\partial D &= B \cdot \Psi^R - \Psi^L \cdot C \\
\partial e_\alpha &= (B\cdot \Psi^R_1 - \Psi^L_2 \cdot C)_{\alpha\alpha}.
\end{align*}
Extend this via the Leibniz rule to obtain a differential on all of $\A$.

\begin{proposition}{\rm\cite{Fr}}\label{prop-invariance}\qua
$\partial^2=0$ and the homology $HC_*(K) = H_*(\A,\partial)$ depends only on the
isotopy class of $K$, not the particular knot diagram. This homology
is called the (framed) knot contact homology of $K$.
\end{proposition}

\noindent There is an equivalence relation on semifree differential
graded algebras known as \textit{stable tame isomorphism}
\cite{Che}, under which $(\A,\partial)$ is independent of the knot
diagram; equivalent DGAs have isomorphic homology. Up to
equivalence, we can refer to $(\A,\partial)$ as the framed knot DGA of
$K$. For a proof of invariance, using a somewhat different
formulation of the framed knot DGA, see \cite{Fr}.

The above definition of the framed knot DGA is combinatorial and
relatively simple, but opaque. A purely topological interpretation
for the full DGA or its homology would be very interesting but is
presently lacking. Most current applications of knot contact
homology use only its lowest-degree component, the degree--$0$
homology $HC_0(K)$, which does have a topological formulation
\cite{KBITwo} as we now describe.

Let $K\subset \R^3 = S^3\setminus\{\text{pt}\}$ be an oriented knot
equipped with the zero framing, and let $l,m$ denote the homotopy
classes of the longitude and meridian of $K$ in ${\pi_1(S^3\setminus
K)}$. Let $\A_K$ denote the tensor algebra over $\Z[\lambda^{\pm
1},\mu^{\pm 1}]$ freely generated by the set ${\pi_1(S^3\setminus
K)}$; a monomial in $\A_K$ looks like $[\gamma_1][\gamma_2]\dots
[\gamma_k]$, where $[\gamma_i]$ denotes the image of
$\gamma_i\in\pi_1(S^3\setminus K)$ in $\A_K$.

\begin{definition}
The
\label{def-htpy}
\textit{cord algebra} of $K$ is the quotient of $\A_K$ by the
relations
\begin{itemize}
\item $[e] = 1+\mu$;
\item $[\gamma m] = [m \gamma] = \mu [\gamma]$ and $[\gamma l] = [l
\gamma] = \lambda [\gamma]$ for $\gamma\in\pi_1(S^3\setminus K)$;
\item $[\gamma_1\gamma_2] + [\gamma_1 m \gamma_2] = [\gamma_1]
[\gamma_2]$ for $\gamma_1,\gamma_2\in\pi_1(S^3\setminus K)$.
\end{itemize}
\end{definition}

If $K$ is the unknot, then $\pi_1(S^3\setminus K) \cong \Z$ is
generated by $m$, and the above relations yield $[e]=1+\mu$ and
$[l]=[e]=\lambda[e]$; it follows easily that the cord algebra of the
unknot is $\Z[\lambda^{\pm 1},\mu^{\pm 1}]/((\lambda-1)(\mu+1))$.

For general knots, the cord algebra is evidently a knot invariant, but
is difficult to compute using \fullref{def-htpy} directly. On the
other hand, degree--$0$ knot contact homology is readily computable in
terms of generators and relations, but is not obviously a topological
invariant. A key result in \cite{Fr} (see also \cite{KBITwo}) states
that these two constructions coincide; we will discuss the intuitive
explanation for this result later in the section.

\begin{proposition}{\rm\cite{KBITwo,Fr}}\label{prop-HCZero}\qua
The cord algebra of $K$ is isomorphic to $HC_0(K)$. Given an
$n$--crossing knot diagram of $K$, the cord algebra can be expressed
as the tensor algebra $\A$ over $\Z[\lambda^{\pm 1},\mu^{\pm 1}]$
generated by $\{a_{ij}\}_{1\leq i,j\leq n,~i\neq j}$ modulo the
$2n^2$ relations given by the entries of the matrices $\Psi^L \cdot
A$ and $A\cdot \Psi^R$.
\end{proposition}

\begin{figure}[ht!]
\centerline{
\includegraphics[height=1.2in]{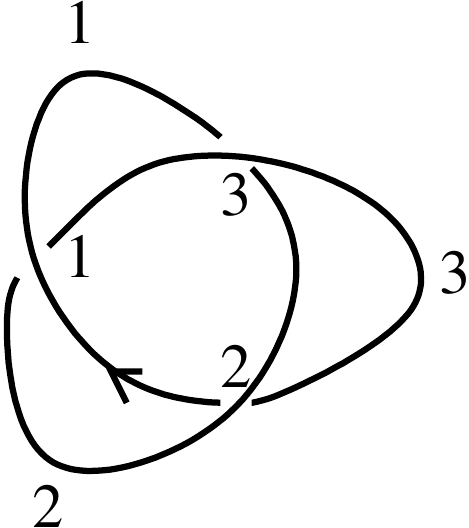}
}
\caption{
The left-handed trefoil, with crossings and diagram components
labeled.
}
\label{fig-LHtrefoil}
\end{figure}

Suppose, for instance, that we wish to compute the cord algebra of
the left-handed trefoil shown in \fullref{fig-LHtrefoil}. We have
$$\Psi^L = \left( \begin{smallmatrix} -a_{21} & \mu & \lambda \\
    1 & -a_{32} & \mu \\ \mu & 1 & -a_{13} \end{smallmatrix}\right),
  \qua\Psi^R = \left( \begin{smallmatrix} -a_{12} & \mu & 1 \\
    1 & -a_{23} & \mu \\ \lambda^{-1}\mu& 1 & -a_{31} \end{smallmatrix}\right),
  \qua\text{and}\ A = \left( \begin{smallmatrix} 1+\mu & a_{12} & a_{13} \\
    a_{21} & 1+\mu & a_{23} \\ a_{31} & a_{32} & 1+\mu
    \end{smallmatrix} \right).$$
When we equate the entries of $\Psi^L\cdot
A$ and $A\cdot\Psi^R$ to zero, we find that $a_{31} = \lambda^{-1}
a_{21}$, $a_{13} = \lambda a_{12}$, $a_{32}=a_{12}$,
$a_{23}=a_{21}$, and $a_{23}=a_{13}$; it follows that we can
substitute for everything in terms of $a_{12}$. If we write
$x=a_{12}$, then the relations reduce to two polynomials in $x$, and
the cord algebra becomes
\[
HC_0(\text{LH trefoil}) \cong \Z[\lambda^{\pm 1},\mu^{\pm 1}][x] /
(\lambda x^2-\lambda x-\mu^2-\mu,\lambda x^2-\mu x-\mu-1).
\]
This is different from the cord algebra of the unknot (set
$\lambda=\mu=1$), and also from the cord algebra of the right hand
trefoil. \fullref{sec-apps} discusses general methods for
distinguishing knot contact homologies and relations to classical
invariants.

The term ``cord algebra'' comes from another interpretation of
$HC_0(K)$. Fix a point $*$ on the knot $K$. A \textit{cord} of $K$
(terminology due to D.\ Bar-Natan) is a path in $S^3$ which begins
and ends on $K\setminus\{*\}$ and does not intersect $K$ in the
interior of the path. Let $\widetilde{\A}_K$ denote the tensor algebra
over $\Z[\lambda^{\pm 1},\mu^{\pm 1}]$ generated by homotopy classes
of cords. In $\widetilde{\A}_K$, mod out by several relations:
\begin{itemize}
\item the contractible homotopy class is equal to $1+\mu$;
\item if two cords are related by pulling an endpoint along
$K$ across $*$, then one is $\lambda$ times the other (the precise
choice is given by orientations);
\item if two cords are related by pushing an interior point in the
cord through $K$, then at the moment when the interior point crosses
$K$, the cord breaks into two, and the four cords involved are
related by
\[
\raisebox{-0.17in}{\includegraphics[width=0.4in]{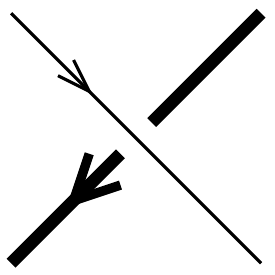}} + \mu
\raisebox{-0.17in}{\includegraphics[width=0.4in]{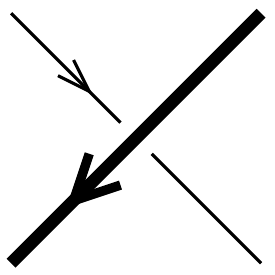}} =
\raisebox{-0.17in}{\includegraphics[width=0.4in]{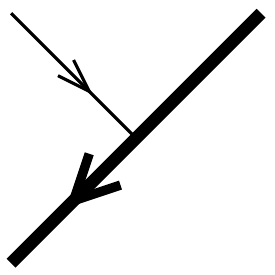}} \cdot
\raisebox{-0.17in}{\includegraphics[width=0.4in]{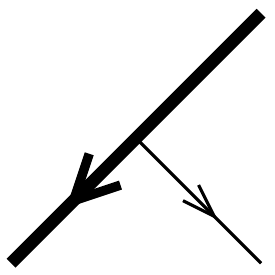}}
\]
where the knot is drawn thickly and the cords thinly.
\end{itemize}

Then $\widetilde{\A}_K$ modulo these relations is isomorphic to the cord
algebra of $K$. To illustrate this isomorphism, push the endpoints
of any cord slightly off of the knot, and join the endpoints via a
path running parallel to the knot, to obtain a closed curve in
$S^3\setminus K$. By carefully tracking the various choices
involved, we can map cords to elements of $\pi_1(S^3\setminus K)$,
and the relations imposed on cords become the relations in
\fullref{def-htpy}.

The isomorphism from \fullref{prop-HCZero} between the cord
algebra and $HC_0(K)$ can now be seen as follows. Place the point
$*$ at the undercrossing of crossing number $1$, and send to
$a_{ij}$ the cord which begins on strand $i$, ends on strand $j$,
and lies strictly above the plane of the diagram in between. The
relations $\Psi^L\cdot A$ and $A\cdot\Psi^R$ which define $HC_0(K)$ are
simply the result of applying the cord ``skein relation'' above when
an endpoint of the cord $a_{ij}$ is pulled through a crossing. See
\cite{Fr} for details.

Cords also give us a heuristic way to see the correspondence between
the cord algebra $HC_0(K)$ and the abstract version of knot contact
homology described in \fullref{sec-LegCH}. If we imagine a cord
as a rubber band and pull it tight, then we can use the cord skein
relation above to write any cord in terms of minimal binormal
chords, that is, line segments which locally minimize distance
between points on the knot. On the other hand, since the DGA which
produces knot contact homology has generators only in nonnegative
degree, any chain in degree $0$ is a cycle; hence the degree--$0$
homology is a quotient of the algebra generated by the degree--$0$
Reeb chords, which are precisely minimal binormal chords.

The results described in this section, and particularly the
interpretation of $HC_0$ as an algebra of cords, can be viewed as a
relative analogue of the work of Viterbo \cite{Vit}, Salamon and
Weber \cite{SW}, and Abbondandolo and Schwarz \cite{AS}, which
shows, roughly speaking, that the symplectic Floer homology of the
cotangent bundle of a manifold is given by the homology of the loop
space of the manifold. In our case, there is a relation between a
relative Floer theory connected to the cotangent bundle and an
algebra generated by ``relative loops'' with endpoints on the knot.
It is likely in this context that knot contact homology is related
to string topology \cite{CS}; see \cite{Coh} for more on the
possible relation in the absolute case.

%*********************************************************************
%*********************************************************************
\section{Knot contact homology: properties}
\label{sec-apps}

Knot contact homology, and in particular the cord algebra $HC_0(K)$,
is a reasonably strong knot invariant.

\begin{proposition}{\rm\cite{KBIOne,Fr}}\qua
As an algebra over $\Z[\lambda^{\pm 1},\mu^{\pm 1}]$, $HC_0(K)$ can
distinguish between mirrors and also between Conway mutants. There
are examples where the cord algebra is stronger than any one of the
following invariants: Alexander polynomial, Jones polynomial, HOMFLY
polynomial, Kauffman polynomial, Khovanov homology,
Ozsv\'ath--Szab\'o invariant.
\label{prop-distinguish}
\end{proposition}

\noindent It is even conceivable at present that the cord algebra
could be a complete knot invariant.

Although the cord algebra is often a difficult object to handle
directly, one can easily compute numerical invariants from it known
as \textit{augmentation numbers}, and
\fullref{prop-distinguish} is proven using these numbers.
Given a ring $R$, which one usually takes to be a finite field
$\Z_p$, augmentation numbers count the number of maps which send
$HC_0(K)$ to $R$, given particular choices of invertible elements in
$R$ for the images of $\lambda$ and $\mu$. The set of possible maps
is finite for given $p$ and can be computed using a program such as
\textit{Mathematica}; the relevant code can be found on the author's
web page.

Knot contact homology has a couple of direct relations to
``classical'' knot invariants.

\begin{proposition}{\rm\cite{Fr}}\qua
One can deduce the Alexander invariant, and thus the Alexander
polynomial, from a linearized form of the framed knot DGA with
$\lambda=1$.
\end{proposition}

\noindent It is also likely that one can find the Alexander
invariant in the cord algebra $HC_0(K)$, without needing to use
higher degree contact homology.

A more interesting link is that between the cord algebra
and the $A$--polynomial of \cite{CCGLS}. The $A$--polynomial is
defined in terms of $SL_2\C$ representations of the knot group,
essentially as follows. For an oriented knot $K$ in $S^3$, the knot
group $\pi_1(S^3\setminus K)$ has two distinguished elements, the
longitude $l$ and meridian $m$ of $K$. Since these two elements
commute, any representation $\rho\co\pi_1(S^3\setminus K) \to
SL_2\C$ is conjugate to one in which $l$ and $m$ are sent to upper
triangular (typically diagonal) matrices
$$\rho(l) = \left(
  \begin{smallmatrix} \lambda & * \\ 0 & \lambda^{-1} \end{smallmatrix}
  \right)\quad\text{and}\quad \rho(m) = \left(
  \begin{smallmatrix} \mu & * \\ 0 & \mu^{-1} \end{smallmatrix}
  \right).$$
The set of pairs of
eigenvalues $(\lambda,\mu)$ over all representations forms a variety
in $(\C^*)^2$ whose $1$--dimensional components are the zero set of
the $A$--polynomial $A_K(\lambda,\mu)$.

One can also form a two-variable polynomial from the cord algebra.
The set of pairs $(\lambda,\mu)$ in $(\C^*)^2$ for which there is a
map $HC_0(K) \otimes \C \to \C$ extending the identity on $\C$ forms
a variety, typically (perhaps always) $1$--dimensional; if it is
$1$--dimensional, then it is the vanishing set of a polynomial
$\widetilde{A}_K(\lambda,\mu)$ which we call the \textit{augmentation
polynomial}. For instance, the left-handed trefoil has augmentation
polynomial given by the resultant of the polynomials $\lambda
x^2-\lambda x-\mu^2-\mu$ and $\lambda x^2-\mu x-\mu-1$ which determine the
cord algebra: up to units, we have
\[
\widetilde{A}_{\text{LH trefoil}}(\lambda,\mu) =
(\lambda-1)(\mu+1)(\lambda-\mu^3).
\]
By contrast, the right-handed trefoil has augmentation polynomial
$$\widetilde{A}_{\text{RH trefoil}}(\lambda,\mu) =
  (\lambda-1)(\mu+1)(1-\lambda \mu^3);$$
this shows that the cord algebra distinguishes mirrors.

The $A$--polynomial is contained in the augmentation polynomial:

\begin{proposition}{\rm\cite{Fr}}\qua
The polynomial
\label{prop-Apoly}
$(1-\mu^2) A_K(\lambda,\mu)$ divides $\widetilde{A}_K(\lambda,-\mu^2)$.
\end{proposition}

\begin{proof}[Sketch of proof]
If we set $\mu=-1$, there is a trivial map from $HC_0(K)$ to $\C$
which sends all $[\gamma]$'s to $0$; it follows that $1-\mu^2$
divides $\widetilde{A}_K(\lambda,-\mu^2)$. Now given a representation
$\rho\co\pi_1(S^3\setminus K) \to SL_2\C$ such that $\rho(l)$ and
$\rho(m)$ are diagonal, define $\{\gamma\} \in \C$, for any
$\gamma\in \pi_1(S^3\setminus K)$, to be the upper left entry of
$\rho(\gamma)$. Then $A_K(\lambda,\mu)=0$, where $\lambda=\{l\}$ and
$\mu=\{m\}$, and
\[
\{A_1A_2\} - \mu \{A_1 \rho(m) A_2\} = (1-\mu^2) \{A_1\} \{A_2\}
\]
for any $A_1,A_2\in SL_2\C$. It follows easily from the homotopy
definition of the cord algebra (\fullref{def-htpy}) that the
map $[\gamma] \mapsto (1-\mu^2)(-\mu)^{\operatorname{lk}(\gamma,K)}
\{\rho(\gamma)\}$ induces a map from $HC_0(K,\lambda,-\mu^2)$ to
$\C$, where $HC_0(K,\lambda,-\mu^2)$ denotes the cord algebra with
$\mu$ replaced by $-\mu^2$. Hence $\widetilde{A}_K(\lambda,-\mu^2)=0$,
so $A_K(\lambda,\mu)$ divides $\widetilde{A}_K(\lambda,-\mu^2)$.
\end{proof}

A recent result of Dunfield and Garoufalidis \cite{DG}, building on
the proof of Property P for knots by Kronheimer and Mrowka
\cite{KM}, states that nontrivial knots have nontrivial
$A$--polynomial. The following result is then a direct consequence of
\fullref{prop-Apoly}.

\begin{corollary}{\rm\cite{Fr}}\qua
Knot contact homology, or in particular the cord algebra,
distinguishes the unknot.
\end{corollary}

We believe that the $A$--polynomial and augmentation polynomial
coincide for two-bridge knots---more precisely, that
$\widetilde{A}_K(\lambda,-\mu^2) = (1-\mu^2) A_K(\lambda,\mu)$ in this
case. This is not, however, true in general. For instance, for the
torus knot $T(3,4)$, we have
\[
\widetilde{A}_{T(3,4)}(\lambda,-\mu^2) = (1-\lambda\mu^8)(1-\mu^2)
A_{T(3,4)}(\lambda,\mu).
\]
It is possible that the augmentation polynomial contains more, or
different, information than the $A$--polynomial.

Since the $A$--polynomial has links to the Jones polynomial (see,
for example, \cite{FGL,Gar}), it is possible that knot contact
homology encodes the Jones polynomial or other classical invariants.
This would make sense given the relation between Chern--Simons knot
invariants and holomorphic curves discussed in the introduction. At
present, however, it is not known how knot contact homology might
contain the Jones polynomial.

\bibliographystyle{gtart}
\bibliography{link}

\end{document}